\documentclass[a4paper]{siamart0216}
\usepackage{subcaption}
\usepackage{booktabs}

\usepackage{geometry}
\usepackage{amsfonts}
\usepackage{graphicx}
\usepackage{epstopdf}
\ifpdf
  \DeclareGraphicsExtensions{.eps,.pdf,.png,.jpg}
\else
  \DeclareGraphicsExtensions{.eps}
\fi

\newcommand{\TheTitle}{An iteration count estimate for a mesh-dependent steepest descent method based on finite elements and Riesz inner product representation} 
\newcommand{\TheAuthors}{Tobias Schwedes, Simon W. Funke, and David A. Ham}

\headers{\TheAuthors}{\TheAuthors}

\title{{\TheTitle}}

\author{
  Tobias Schwedes\thanks{Department of Mathematics, Imperial College London, London, UK and Grantham Institute, Imperial College London, London, UK (\email{t.schwedes14@imperial.ac.uk}).}
  \and
  Simon W. Funke\thanks{Center for Biomedical Computing, Simula Research Laboratory, Oslo, Norway 
    (\email{simon@simula.no}).}
      \and
  David A. Ham\thanks{Department of Mathematics, Imperial College London, London, UK (\email{david.ham@imperial.ac.uk}).}
}

\usepackage{amsopn}

\ifpdf
\hypersetup{
  pdftitle={\TheTitle},
  pdfauthor={\TheAuthors}
}
\fi

\begin{document}

\maketitle

\begin{abstract}
Existing implementations of gradient-based optimisation methods typically
assume that the problem is posed in Euclidean space. When solving
optimality problems on function spaces, the functional derivative is then inaccurately represented with respect to $\ell^2$ instead of the inner product induced by the function space. This error manifests as a mesh dependence in the number of iterations required to solve the optimisation problem. In this paper, an analytic estimate is derived for this iteration count in the case of a simple and generic discretised optimisation problem. The system analysed is the steepest descent method applied to a finite element problem.  The estimate is based on Kantorovich's inequality and on an upper bound for the condition number of Galerkin mass matrices. Computer simulations validate the iteration number estimate. Similar numerical results are found for a more complex optimisation problem constrained by a partial differential equation. Representing the functional derivative with respect to the inner product induced by the continuous control space leads to mesh independent convergence.
\end{abstract}



\section{Introduction}

Computationally solving an optimisation problem on Hilbert spaces requires
discretisation and the application of an optimisation method to the resulting
finite dimensional problem. In this context, mesh independent convergence of
the optimisation algorithm means that, for a discretisation given by a
sufficiently fine mesh, the number of iterations required to solve the
optimisation problem to a given tolerance is bounded. Conversely, for a mesh
dependent algorithm there exist arbitrarily fine meshes which will require an
arbitrarily large number of iterations to achieve convergence.

As we will discuss in section \ref{sec:related}, employing a ``suitable'' inner
product in the optimisation algorithm is a necessary condition for achieving
mesh independent convergence. However, many well-established continuous optimisation packages are
not designed for function based optimisation and are hard-coded to apply the Euclidian $(\ell^2)$ inner product
(for example the MATLAB Optimization Toolbox \cite{matlab2016}, TAO \cite{benson2005tao}, SNOPT \cite{gill2005snopt}, scipy.optimize \cite{scipy2001,byrd1995limited} and IPOPT \cite{kawajir2015,wachter2006implementation}). One can intuitively comprehend the inaccuracy in this practise for function based optimisation: since the inner product defines angles and distances, drawbacks in the convergence of optimisation methods using an inner product associated with a space different from the control space are expected. For example, this error manifests in suboptimal search directions. For gradient based optimisation methods, drawbacks can similarly be expected from the fact that the gradient representation of the derivative associated with the functional of interest is inner product dependent. As a result, one obtains mesh
dependent convergence rates when these packages are applied to problems on
Hilbert spaces. Only recently inner product aware optimisation packages have appeared \cite{young2014, funke2016}. Thus, given the prevalence of this practise, it is germane to attempt to quantify the
extent of mesh dependence in optimisation problems. In other words, to what
extent is mesh independence a mathematical nicety, and to what extent an
essential component of an acceptable optimisation algorithm? To the best of our
knowledge, there are no analytical results on the iteration count of mesh
dependent iterative methods in existing literature.\\

In this paper, we study the finite element discretisation of a generic
optimisation problem formulated in a Hilbert function space. If
we solve this problem using the steepest descent method
with respect to the $\ell^2$ inner product, the
convergence rate is mesh dependent. For this situation we analytically derive an
iteration number estimate that reveals that the number of iterations is polynomial in the ratio between largest and smallest
directional element size of the mesh. Conversely, if the optimisation
is conducted in the inner product induced by the Hilbert space,
the optimisation converges in exactly one iteration, independently of
the computational mesh employed. Numerical experiments confirm the
analytical estimate in the special case of the Hilbert space $L^2$.

The estimate is derived for a particularly simple case in the
expectation that the impact of an ``incorrect'' inner
product on more complex problems will be at least as severe. We illustrate this
using the example of the optimal control of the Poisson equation, a partial differential equation (PDE)
constrained optimisation problem. When the $\ell^2$ inner product is employed in the
optimisation, the iteration count is mesh dependent, similar to that observed
in the simple case. Conversely, by using an implementation  that employs the
$L^2$ inner product, mesh independent convergence is achieved.

This work is organised as follows: In section \ref{sec:related}, we examine the existing work on mesh dependence. Section \ref{sec:prelim} restates the Riesz representation theorem from functional analysis. In section \ref{sec:formulation}, the continuous optimisation problem considered throughout this paper is formulated, and section \ref{sec:finite_element} is dedicated to the finite element discretisation of the problem. In section \ref{sec:iteration_L2}, the iteration count for the steepest descent making use of the $L^2$ represented Riesz gradient is derived, while  section \ref{sec:iteration_l2} develops an estimate for the iteration number when a gradient is used that is represented with respect to the Euclidean scalar product. Numerical experiments investigating the validity in applications for the stated iteration counts are presented in section \ref{sec:experiments}, as well as a simulation study considering mesh dependence for a PDE constrained problem.


\section{Related work}\label{sec:related}

Recent work on mesh independence in finite element
methods focusses on achieving this characteristic in the solution of
the primary finite element problem
\cite{kirby2010functional,gunnel2014note,mardal2011preconditioning}. This
is, in fact, very closely related to the problem considered here:
solving the finite element problem amounts to minimising the residual,
and the task of mapping the residual gradient from the dual space back
into the primal function space is similar to employing the correct
inner products in calculating updates for the optimisation problem.

\cite{kirby2010functional} derives bounds for the condition number of
finite element operator matrices based on the structure that the
discretisation inherits from the underlying continuous Hilbert
space. The condition number depends on the underlying discretisation
due to a non-isometric embedding of the operator into Euclidean
space. Spectral bounds that are independent of the finite element
subspace (that is in particular of the mesh) are achieved by applying the Riesz map as
preconditioner. This abstract preconditioning principle transfers
directly to the problem considered here: mesh independent convergence
requires that the Fr\'echet derivative of the objective functional,
defining the search direction in the optimisation, is represented with
respect to the correct Hilbert space. This is equivalent to employing
an $\ell^2$ representation and preconditioning with the Riesz map. In
the case of the Hilbert space $L^2$, the Riesz map of the discretised
problem is the inverse of the Galerkin mass matrix.

More generally, \cite{gunnel2014note} shows that the choice of a preconditioner is equivalent to the choice of an inner product on the underlying space. That work investigates the dependence of conjugate gradient and MINRES methods in Hilbert spaces on the underlying scalar product. In finite dimensions, the naturally induced scalar product is related to every other inner product by the application of a positive definite matrix. The Riesz map preconditioner corresponds to an application of that matrix's inverse. For $L^2$, this is once again inverse of the Galerkin mass matrix.\\

 \cite{kirby2010functional,gunnel2014note,mardal2011preconditioning} state conditions under which mesh dependence can be avoided rather than discussing the drawbacks of the opposite. Considering that many well-established optimisation packages assume that the problem is formulated in Euclidean space, we instead investigate of the magnitude of the resulting performance loss.

\section{Preliminaries}\label{sec:prelim}

Let $(H, \langle \cdot , \cdot {\rangle}_H)$ be a Hilbert space and let $f:H\rightarrow \mathbb{R}$ be Fr\' echet differentiable. For any $u\in H$, let $f'(u)\in \mathcal{L}(H, \mathbb{R})=H^*$ be the Fr\' echet derivative of $f$ at $u$. According to the Riesz representation theorem, there is a unique $\hat{u} \in H$ such that
\begin{align}
 f'(u) v  = \langle \hat{u}, v {\rangle}_{H}
\label{eq_riesz-rep-thm}
\end{align}
for all $v \in H$. The map $\mathcal{R}_H: H^* \rightarrow H$ that sends $f'(u)$ to its unique representative $\hat{u}\in H$ is called the Riesz map. Further, $\mathcal{R}_{H}(f'(u)))=\hat{u}$ is called the Riesz representer of $f'(u)$ in $H$, or the gradient of $f$ in $u$ represented with respect to $H$.

\section{Formulation}\label{sec:formulation}

In order to discuss how the choice of the inner product for the representation of the functional derivative impacts upon the convergence of the optimisation with respect to the underlying discretisation, a problem as simple as possible, at the same time being generic in this context, is considered. Such an optimisation problem is then exemplary while its simplicity suggests dependence on the inner product and the discretisation of at least similar order for more complex optimisation problems, including problems with PDE constraints. Consider the minimisation problem,
\begin{align}
\min_{u \in H} \hspace{1mm} \alpha \langle u,u {\rangle}_H +\beta \langle u, v {\rangle}_H + \gamma,
\label{abstract_min_problem}
\end{align}
where $\alpha>0$, $\beta, \gamma \in \mathbb{R}$ and $v \in H$. Further, $H\subset \{u: \Omega \rightarrow \mathbb{R}\}$, with bounded domain $\Omega\subset \mathbb{R}^n $ and $n \in \mathbb{N}$, is a Hilbert function space with inner product $\langle \cdot, \cdot \rangle_H$ that is of some integral type. For instance, $H$ could be a Hilbert Sobolev space, i.e$.$ $H=W^{k,2}$ with $k\in\mathbb{N}$. The derivation of the steepest descent iteration count estimate for \eqref{abstract_min_problem} can be reduced to the consideration of the simple optimisation problem
\begin{align}
\min_{u\in L^2(\Omega)}  \left\{ f(u) = \frac{1}{2} \langle 1-u, 1-u {\rangle}_{L^2} \right\},
\label{example_opti}
\end{align}
where $\langle \cdot, \cdot {\rangle}_{L^2}$ denotes the inner product of the space $L^2:=L^2(\Omega)$. This simplification is justified by the common finite element approach and the quadratic structure of both discretisations. Details are found in section $5$. It is easy to see that \eqref{example_opti} is well posed with $u=1$ as unique solution. \\


In the subsequent sections we consider the solution of \eqref{example_opti} by using the steepest descent algorithm with exact line search and the finite element method. In doing so, two approaches are compared: In the first one, the continuous formulation of $f$ in \eqref{example_opti} is employed in order to compute the gradient used for the optimisation method. Then, the finite element discretisation is applied. In the second one, the discretisation is performed first, whereupon the gradient is computed with respect to the $l^2$ inner product of the coefficient vectors. We will see that the formulation used in the second approach is mesh dependent by a scaling of the local finite element mass matrix. Results of the optimisation procedure using both approaches are compared for non-uniform meshes.




\section{Finite element discretisation}\label{sec:finite_element}

Let us now consider the nodal finite element discretisation of the continuous optimisation problem \eqref{example_opti}. Let $\mathcal{T}$ be a tessellation of the domain $\Omega$ by topologically identical $n$-polytopes $K$. Typically the cells, $K$ are simplices or hypercubes. Assume that any $K\in \mathcal{T}$ is equipped with an individual set of linearly independent functions $\mathcal{P}=\{\phi_1,...,\phi_p\}\subset L^2(K)$ and another set of functions $\mathcal{N}=\{N_1,...,N_p\}$, where $N_1,...,N_p:L^2(K)\rightarrow \mathbb{R}$ are called nodal variables and form a basis of $\mathcal{P}'$, the dual space of $\mathcal{P}$. Then, each $(K, \mathcal{P}, \mathcal{N})$ forms a finite element following Ciarlet's definition in \cite{ciarlet2002finite}. If $v$ is a function in $L^2(K)$,
the local interpolant of $v$ on $K$ is defined by
\begin{align}
\mathcal{I}_K(v) = \sum_{i=1}^p N_i(v) \phi_i,
\label{local_interpolant}
\end{align}
following the formulation in \cite{brenner08}. Thus, $v$ is approximated by the sum of basis functions $\phi_1,...,\phi_p$ and coefficients $N_1(v),$ $..., N_p(v)$. Typically, the application of $N_1,...,N_p$ to $v$ involves the evaluation of $v$ or one of its derivatives on some given point on $K$. Let $u\in V \subset C^m(\Omega)\subset L^2$, where $m$ is the order of the highest partial derivatives involved in evaluating $N_i(\left.u\right|_{K})$ for $i=1,...,p$. Then, the global interpolant of $u$ on $\mathcal{T}$ is defined by
\begin{align}
\left.\mathcal{I}_{\mathcal{T}}(u) \right|_{K} = \mathcal{I}_K(\left. u\right|_K)
\end{align}
for any $K \in \mathcal{T}$. Note that $\mathcal{I}_{\mathcal{T}}$ may have multiple values at the interfaces of elements. However, this is not problematic since these sets have Lebesgue measure equal to zero. It is obvious that $\mathcal{I}_{\mathcal{T}}(u)$ can be written similarly to \eqref{local_interpolant} when $\phi_1,...,\phi_p$ are expanded to $\Omega$ as follows: in case $\left.\phi_i \right|_{\partial K}=0$, set $\phi_i$ equal to zero outside of $K$. If $\phi(x)\neq 0$ for some $x\in\partial K$, set $\phi_i$ equal to the basis functions on neighbouring elements that are also not equal to zero in $x$. On all other elements, set $\phi_i$ equal to zero.
The extensions of the basis functions of all elements are given by the so called global basis functions $\varphi_1,...,\varphi_d$, and their corresponding global nodal variables by $N^{\varphi}_1,...,N^{\varphi}_d$, where
\begin{align}
u_i := N^{\varphi}_{i} (u) = N_{\hat{i}}(\left. u\right|_K),
\end{align}
for $\left.\varphi_{i}\right|_K = \phi_{\hat{i}} \neq 0$ on $K$, where $\phi_{\hat{i}}$ denotes the basis function associated with the global basis function $\varphi_i$ on $K$. Hence, the global interpolant $\mathcal{I}_{\mathcal{T}}(u)\in V_h \subset V$ with $\dim(V_h)=d$ can be expressed as
\begin{align}
\mathcal{I}_{\mathcal{T}}(u)(x) = \sum_{i=1}^d u_i \varphi_i(x).
\label{global_interpolant}
\end{align}
One computes 
\begin{align}
f({\mathcal{I}_{\mathcal{T}}(u)}) &=  \int_{\Omega} \frac{1}{2}\Big(1_i-\mathcal{I}_{\mathcal{T}}(u)(x)\Big)^2 \mathrm{d}x \\
&= \frac{1}{2} \int_{\Omega} \Big(\sum_{i=1}^d (1_i-u_i) \varphi_i(x)\Big)^2 \mathrm{d}x \\
&= \frac{1}{2} \sum_{i,j=1}^d (1_i-u_i) (1_i-u_j) \int_{\Omega} \varphi_i(x) \varphi_j(x) \mathrm{d}x,
\end{align}
where $1_i = N_i^{\varphi}(1)$ for any $i=1,...,d$. Given the mass matrix $M$ with $M_{ij} := \int_{\Omega} \varphi_i \varphi_j \mathrm{d}x$, the finite element discretised version of the original optimality problem \eqref{example_opti} can be expressed as
\begin{align}
\min_{\mathcal{I}_{\mathcal{T}}(u)} f(\mathcal{I}_{\mathcal{T}}(u)) =  \min_{\vec{u}\in \mathbb{R}^d} \Big\{ f(\vec{u}) = \frac{1}{2} ( \vec{1} - \vec{u})^T M (\vec{1} - \vec{u})\Big\},
\label{discretised_opti}
\end{align}
where $\vec{u}= (u_1,\ldots, u_d), \vec{1} = (1_1,\ldots, 1_d) \in \mathbb{R}^d$. Here, we identified $\left.f \right|_{V_h}$ with the function
\begin{align}
\vec{u} \mapsto \frac{1}{2} ( \vec{1} - \vec{u})^T M (\vec{1} - \vec{u}).
\end{align}
This makes sense since for any $u\in L^2$ there is a unique $\vec{u}\in \mathbb{R}^d$ such that \eqref{global_interpolant} holds true.

\section{Iteration count using $L^2$ inner product}\label{sec:iteration_L2}

In order to solve problem \eqref{example_opti}, we apply the steepest descent method while representing the gradient according to Riesz theorem with respect to the $L^2$ inner product. Initialized by $u_0$, the first iterate $u_1$ is given by $u_1 = u_0 - \alpha \mathcal{R}_{L^2}(f'(u_0))$, where $\mathcal{R}_{L^2}(f'(u_0))$ denotes the Riesz representer of the first Fr\' echet derivative $f'(u_0)$ in $L^2$. We compute $f'(u)\in \mathcal{L}(L^2, \mathbb{R})$ as
\begin{align}
f'(u)(\cdot) &= -\frac{1}{2} \Big(\langle \cdot, 1-u {\rangle}_{L^2} + \langle 1-u, \cdot {\rangle}_{L^2}\Big) \\
&= - \langle 1-u, \cdot {\rangle}_{L^2},
\end{align}
where we used the symmetry of $\langle\cdot, \cdot{\rangle}_{L^2}$. Hence,
\begin{align}
\mathcal{R}_{L^2}\left(f'(u)\right) = -(1-u).
\label{L_2Rieszrep}
\end{align}
The second Fr\' echet derivative $f''(u) \in \mathcal{L}(L^2, \mathcal{L}(L^2, \mathbb{R}))$ is given by
\begin{align}
f''(u) \cdot \cdot = \langle \cdot, \cdot {\rangle}_{L^2}.
\end{align}
For any $n>2$, we find $F^{(n)}(u) = 0$. Applying Taylor's theorem for function spaces at $u_1 = u_0 - \alpha \mathcal{R}(f'(u_0))$ using \eqref{L_2Rieszrep} yields
\begin{align}
f(u_1) &= f(u_0) - \alpha f'(u_0)\mathcal{R}_{L^2}\left(f'(u_0)\right) + \frac{1}{2}\alpha^2 f''(u_0)\left(\mathcal{R}_{L^2}\left(f'(u_0)\right)\right)^2 \\
&= \left(\frac{1}{2} - \alpha + \frac{1}{2}\alpha^2\right)\|1-u_0\|^2_{L^2}.
\end{align}
One easily computes that the minimum of $f(u_1)$ is found for $\alpha^* = 1$, i.e$.$ $\alpha^*$ is the step size of an exact line search, and $u_1 \equiv 1$. Hence, the solution of the optimality problem \eqref{example_opti} is found after one steepest descent iteration, independently of the initial $u_0$. This still holds for the finite element discretised version \eqref{discretised_opti} of problem \eqref{example_opti}. This can be seen from the fact that
\begin{align}
\mathcal{R}_{L^2}\left(f'(u_0)\right) = u_0-1 = \sum_{i=1}^d (u_0-1)_i \cdot \varphi_i = \mathcal{I}_{\mathcal{T}}\left(\mathcal{R}_{L^2}\left(f'(u_0)\right)\right),
\end{align}
for any $u_0 = \sum_{i=1}^d (u_0)_i \varphi_i \in \mathcal{I}_{\mathcal{T}}(L^2):= \{\mathcal{I}_{\mathcal{T}}(u): u \in L^2\}$, where $(u_0-1)_i:= (u_0)_i - 1_i$ for any $i=1,...,d$. Hence, the first steepest descent iterate for the discretised optimisation problem \eqref{discretised_opti} equals the constant function $1$, which is the optimum.

\section{Iteration count estimate using $\ell^2$ inner product}\label{sec:iteration_l2}

In this section, an iteration count estimate is derived for the steepest descent method using gradients with respect to the $\ell^2$ inner product,  computed based on the finite element discretised version \eqref{discretised_opti} of the optimality problem \eqref{example_opti}. Since $M$ is symmetric, one computes
\begin{align}
f'(\vec{u}) \cdot = -(\vec{1}-\vec{u})^T M \cdot &= \langle - M(\vec{1}-\vec{u}), \cdot {\rangle}_{\ell^2},
\end{align}
where $\langle\cdot, \cdot{\rangle}_{\ell^2}$ denotes the $\ell^2$ inner product, i.e$.$ $\langle \vec{u}, \vec{v} {\rangle}_{\ell^2}:= \vec{u}^T v = \sum_{i=1}^d u_i v_i$ for any $\vec{u}, \vec{v} \in \mathbb{R}^d$. Thus,
\begin{align}
\mathcal{R}_{\ell^2}(f'(\vec{u})) = \nabla f(\vec{u}) = -M(\vec{1}-\vec{u}).
\label{littlel2_gradient}
\end{align}
Notice that $\mathcal{R}_{\ell^2}(f')$ differs from $\mathcal{R}_{L^2}(f')$ in \eqref{L_2Rieszrep} by the multiplication of the mass matrix $M$. Since $M$ reflects the structure of the underlying mesh, that is, the spatial distribution of elements and their sizes, it is expected that the convergence of the steepest descent method using $\mathcal{R}_{\ell^2}(f')$ is mesh dependent.

Using \eqref{littlel2_gradient}, the $(k+1)$th iterate of the steepest descent method is given by $\vec{u}_{k+1} = \vec{u}_k + \alpha_k M(\vec{1}-\vec{u}_k)$. In order to minimise
\begin{align}
f(\vec{u}_{k+1}) = \frac{1}{2} \Big(\vec{1} - (\vec{u}_k + \alpha_k M(\vec{1}-\vec{u}_k) \Big)^T M\Big(\vec{1} - (\vec{u}_k + \alpha_k M(\vec{1}-\vec{u}_k)\Big),
\end{align}
with respect to $\alpha_k$, we set $\frac{\partial f}{\partial \alpha_k}(\vec{u}_{k+1})=0$ such that
\begin{align}
\alpha_k &=  \frac{(\vec{1}-\vec{u}_k)^T M^2 (\vec{1}-\vec{u}_k)}{(\vec{1}-\vec{u}_k)^T M^3 (\vec{1}-\vec{u}_k)}\label{eq:alpha_exact}
\end{align}
Note that the Hessian of $f$ equals $M$. Considering the second derivative of $f(\vec{u}_{k+1})$ with respect to $\alpha_k$ and using the fact that mass matrices are positive definite shows that $\alpha_k$ as 	given in \eqref{eq:alpha_exact} is the unique minimiser, i.e$.$ $\alpha_k$ is the step size of an exact line search.\\

Using steepest descent with exact line search on a strongly convex 
quadratic function, we may apply the Lemma of Kantorovich from section 8.6 in \cite{luenberger08}. The lemma provides an recursive error estimate in $f$ at the $k$th iterate with respect to the condition number $\kappa(M)$ of $M$. Since $M$ is normal, its condition number is given by
\begin{align}
\kappa(M) = \frac{\lambda_{\max}^M}{\lambda_{\min}^M},
\end{align}
where $\lambda_{\max}^M$ and $\lambda_{\min}^M$ denote the maximum and minimum eigenvalues of $M$, respectively. Applying the Lemma of Kantorovich yields
\begin{align}
f(\vec{u}_{k}) - f(\vec{u}^*) &\le \left[\frac{\kappa(M)-1}{\kappa(M)+1}\right]^2 \Big(f(\vec{u}_{k-1}) - f(\vec{u}^*)\Big) \label{1st_inequality_kantorovich} \\
&\le \ldots \nonumber\\
&\le \left[\frac{\kappa(M)-1}{\kappa(M)+1}\right]^{2k} \Big(f(\vec{u}_0) - f(\vec{u}^*)\Big)  \\
&= \left[\frac{\kappa(M)-1}{\kappa(M)+1}\right]^{2k} f(\vec{u}_0), \label{eq:error_est_steep}
\end{align}
where $\vec{u}^*=(1_1,...,1_d)\in \mathbb{R}^d$ is the optimal value such that $f(\vec{u}^*)=0$. \\

An obvious question is how to determine $\kappa(M)$ from the underlying discretisation of the domain. Corollary 1 of \cite{fried72} gives an upper bound for the condition number of mass matrices depending on the minimum and maximum eigenvalues of their local mass matrices. Let $M_K$ be the local mass matrix associated with element $K$, and let $\lambda_{\max}^{M_K}$ and $ \lambda_{\min}^{M_K}$ denote the maximum and minimum eigenvalues of $M_K$, respectively. Applying the estimate to $M$, we obtain
\begin{align}
\kappa(M) \le p_{\max} \frac{\max_{K} \lambda_{\max}^{M_K}}{\min_{K} \lambda_{\min}^{M_K}},
\label{eq:condition_estimate}
\end{align}
where $p_{\max}$ is the maximum number of elements around any nodal point. In order to determine $\lambda_{\max}^{M_K}$ and $\lambda_{\min}^{M_K}$, respectively, we consider the bijective transformation $T_K$ from a reference element $(\hat{K}, \hat{\mathcal{P}}=\{\psi_1, ..., \psi_p\},$ $\hat{\mathcal{N}}=\{\hat{N}_1,...,\hat{N}_p\})$ to the local element $(K,\mathcal{P}^K=\{\phi_1^K, ..., \phi_p^K\},\mathcal{N}^K=\{N_1^K,...,N_p^K\})$, where $\phi_l^K\circ T_K = \psi_l $. Let us assume affine transformations, i.e$.$
\begin{equation}
\begin{aligned}
T_K: \hat{K} &\rightarrow K \\
x &\mapsto J_K x + y_K,
\end{aligned}
\end{equation}
for $J_K\in\mathbb{R}^{n\times n}$ invertible and $y_K\in \mathbb{R}^n$. Applying singular value decomposition to the Jacobian matrix $J_K$ of the transformation $T_K$, there are unitary $U, V \in \mathbb{R}^{n \times n}$ such that
\begin{align}
J_K = U H_K V,
\end{align}
where $H_K\in\mathbb{R}^{n\times n}$ is a diagonal matrix with positive real numbers $h_1^K,...,h_n^K$ on the diagonal. Thus,
\begin{align}
|\det(J_K)| = |\det(H_K)|  = h_1^K \cdot ... \cdot h_n^K.
\end{align}
The factors $h_1^K,...,h_n^K$ describe the scaling of lengths in orthogonal directions between reference element $\hat{K}$ and local element $K$, respectively. Using the transformation theorem, an entry $\mu_{lm}^K$ with $l,m=1,...,p$ of the local mass matrix $M_K$ associated with element $K$ can therefore be written as
\begin{align}
\mu_{lm}^K&= \int_K \phi_l^K(x) \phi_m^K(x) \mathrm{d}x \label{element_transformation_1}\\
&= \int_{\tilde{K}} \phi_l^K(T_K(\zeta)) \phi_m^K(T_K(\zeta)) |\det(J_K)| \mathrm{d}\zeta \label{element_transformation_2}\\
&= h_1^K \cdot ... \cdot h_n^K \int_{\tilde{K}} \psi_{l}(\zeta) \psi_{m}(\zeta) \mathrm{d}\zeta \label{element_transformation_3}\\
&= h_1^K \cdot ... \cdot h_n^K \cdot \hat{\mu}_{lm}.
\label{element_transformation_4}
\end{align}
Consequently, the local mass matrix $M_K$ has the form
\begin{align}
M_K = h_1^K \cdot ... \cdot h_n^K  \hat{M},
\label{prop_matrices}
\end{align}
where $\hat{M} = ( \hat{\mu}_{lm} )_{l, m=1,...,p}$ is the mass matrix associated with the reference element $\hat{K}$. Hence, the eigenvalues of $M_K$ are proportional to the eigenvalues of $\hat{M}$, with proportionality factor equal to $h_1 \cdot ... \cdot h_n$. Note that the matrix $\hat{M}$ is independent of the underlying element $K$. \\


Let $\lambda_{\max}^{\hat{M}}$ and $\lambda_{\min}^{\hat{M}}$ denote the largest and smallest eigenvalues of $\hat{M}$, respectively. The estimate \eqref{eq:condition_estimate} can now be written as
\begin{align}
\kappa(M) &\le p_{\max} \frac{\lambda_{\max}^{\hat{M}}}{\lambda_{\min}^{\hat{M}}} \frac{\max\left\{\prod_{i=1}^n h_i^K: K\right\}}{\min\left\{\prod_{i=1}^n h_i^K: K\right\}} \\
&\le p_{\max} \frac{\lambda_{\max}^{\hat{M}}}{\lambda_{\min}^{\hat{M}}} \prod_{i=1}^n \frac{\max_K h_i^K}{\min_K h_i^K}.
\label{new_estimate}
\end{align}
Due to the role of $h_1^K, ..., h_p^K$ as scaling factors of the element $K$ relative to the reference element $\tilde{K}$ in orthogonal directions, the expression $\max\left\{\prod_{i=1}^n h_i^K: K\right\}/ \min\left\{\prod_{i=1}^n h_i^K: K\right\}$ may be considered as a measure of non-uniformity in the mesh.\\

In the particular case where the maximum and minimum of the scaling factors $h_i^K$ in orthogonal directions coincide for $i=1,...,n$, we may set $h_{\max}:=\max_K h_i^K$ and $h_{\min}:=\min_K h_i^K $ for any $i$, respectively, such that \eqref{new_estimate} simplifies to
\begin{align}
\kappa(M) &\le p_{\max}\frac{\lambda_{\max}^{\hat{M}}}{\lambda_{\min}^{\hat{M}}} \left(\frac{h_{\max}}{h_{\min}}\right)^n.
\end{align}
Assuming that the underlying mesh contains significant non-uniformity, i.e$.$ $h_{\max}/h_{\min} \gg 1$, the following approximation holds true:
\begin{align}
\log \left( \frac{p_{\max}\frac{\lambda_{\max}^{\hat{M}}}{\lambda_{\min}^{\hat{M}}} \left(\frac{h_{\max}}{h_{\min}}\right)^n-1}{p_{\max}\frac{\lambda_{\max}^{\hat{M}}}{\lambda_{\min}^{\hat{M}}} \left(\frac{h_{\max}}{h_{\min}}\right)^n+1} \right) &= \sum_{\nu=1}^{\infty} \frac{(-1)^{\nu}}{\nu}\left(\frac{p_{\max}\frac{\lambda_{\max}^{\hat{M}}}{\lambda_{\min}^{\hat{M}}} \left(\frac{h_{\max}}{h_{\min}}\right)^n-1}{p_{\max}\frac{\lambda_{\max}^{\hat{M}}}{\lambda_{\min}^{\hat{M}}} \left(\frac{h_{\max}}{h_{\min}}\right)^n+1}-1\right)^{\nu}\\
&= \sum_{\nu=1}^{\infty} - \frac{2^{\nu}}{\nu}\cdot\left(p_{\max}\frac{\lambda_{\max}^{\hat{M}}}{\lambda_{\min}^{\hat{M}}} \left(\frac{h_{\max}}{h_{\min}}\right)^n+1\right)^{-\nu} \\
&\approx -2 \cdot \left(p_{\max}\frac{\lambda_{\max}^{\hat{M}}}{\lambda_{\min}^{\hat{M}}} \left(\frac{h_{\max}}{h_{\min}}\right)^n+1\right)^{-1}.
\end{align}
Setting \eqref{eq:error_est_steep} smaller or equal to $\varepsilon$, one then obtains
\begin{align}
k \ge -\frac{1}{4} \log\left( \varepsilon/f(\vec{u}_0)\right) \cdot \left( p_{\max}\frac{\lambda_{\max}^{\hat{M}}}{\lambda_{\min}^{\hat{M}}} \left(\frac{h_{\max}}{h_{\min}}\right)^n +1\right) =: \hat{k}.
\label{final_iters_estimate}
\end{align}
Thus, after $\hat{k}$ iterations the error of the steepest descent method is smaller or equal to $\varepsilon$. Careful interpretation of this result yields that an iteration number satisfying that the associated error is smaller or equal to $\varepsilon$ is at most polynomial in the ratio between largest and smallest directional element size of order equal to the domain dimension $n$. The subsequent simulations show that this dependency is also observed numerically.\\

\textbf{Remark:} In what follows, we briefly justify the simplification of considering \eqref{abstract_min_problem} instead of \eqref{example_opti} in the context of deriving \eqref{final_iters_estimate}. With the notation used above, discretising
\begin{align}
\alpha \langle u,u {\rangle}_H +\beta \langle u, v {\rangle}_H + \gamma,
\end{align}
using finite elements leads to
\begin{align}
\alpha \sum_{i,j=1}^d u_i u_j \langle \varphi_i, \varphi_j {\rangle}_H  + \beta \sum_{i,j=1}^d u_i v_j \langle \varphi_i, \varphi_j {\rangle}_H + \gamma = \alpha \vec{u}^T M_H \vec{u} + \beta \vec{u}^T M_H \vec{v} +\gamma,
\label{fem_discr_abstract_problem}
\end{align}
where $v_1, ...,v_n$ are the coefficients of the interpolant $\mathcal{I}_{\mathcal{T}}(v) = \sum_{i=1}^n v_i \varphi_i$ of $v\in H$, and $M_H = (\langle \varphi_i, \varphi_j {\rangle}_H)_{i,j}$. In the special case of $H=L^2([0,1])$, it holds that $M_H$ coincides with the Galerkin mass matrix $M$. Since adding constants does not change the solution of an optimisation problem, minimising the right hand side of \eqref{fem_discr_abstract_problem} is equivalent to minimising
\begin{align}
\frac{1}{2}\vec{u}^T \tilde{M}_H \vec{u} + \vec{u}^T \vec{\tilde{v}},
\label{calc1_derivation_simplification}
\end{align}
where $\tilde{M}_H = 2\alpha M_H$ and $\vec{\tilde{v}}= \beta M_H \vec{v}$. Analogously, minimising \eqref{calc1_derivation_simplification} is equivalent to minimising
\begin{align}
 \frac{1}{2}\vec{u}^T \tilde{M}_H \vec{u} + \vec{u}^T \vec{\tilde{v}} + \frac{1}{2}{\vec{u}^*}\hspace{0mm}^T \tilde{M}_H \vec{u}^* = \frac{1}{2} (\vec{u} - \vec{u}^*)^T \tilde{M}_H (\vec{u}- \vec{u}^*),
\label{calc2_derivation_simplification}
\end{align}
where $\vec{u}^*$ is the optimal solution, for which obviously $\tilde{M}_H \vec{u}^*= -\vec{\tilde{v}}$. Since Kantorovich's inequality is valid for any $\vec{u}^*$ and the estimate in \eqref{1st_inequality_kantorovich} depends only on the condition number of the Hessian matrix of the quadratic problem, one may conclude that the iteration count estimate is valid for minimising \eqref{calc2_derivation_simplification} in the same way as for the minimisation of
\begin{align}
\frac{1}{2} (\vec{u} - \vec{1})^T \tilde{M}_H (\vec{u} - \vec{1}).
\label{new_simple_problem}
\end{align}
We need to show that \eqref{final_iters_estimate} applies for $\tilde{M}_H$ similarly as for $M$. For the sake of clarity, substitute $\tilde{M}_H$ with $M_H$. Let $\langle \cdot, \cdot {\rangle}_{H(K)}$ and $\langle \cdot, \cdot {\rangle}_{H(\hat{K})}$ be defined as the inner products of the same integral type as $H$, with integration domains equal to $K$ and $\hat{K}$, respectively. The local finite element matrix ${M}_{H,K}$ corresponding to ${M}_H$ is then given by the entries
\begin{align}
\mu_{lm}^{H,K}= 2\alpha \cdot \langle \phi_l^K, \phi_m^K {\rangle}_{H(K)},
\end{align}
where $\phi_l^K$ with $l=1,...,p$ denote the local basis functions
associated with element $K$. In a manner analogous to the proof of Theorem 1 in \cite{fried72}, one can show that the condition number estimate \eqref{eq:condition_estimate} holds true for ${M}_H$, i.e$.$
\begin{align}
\kappa({M}_H) \le p_{\max} \frac{\max_K \lambda_{\max}^{{M}_{H,K}}}{\min_K \lambda_{\min}^{{M}_{H,K}}}.
\end{align}
The reference finite element matrix $\hat{M}_H$ corresponding to ${M}_H$ is given by the entries
\begin{align}
\hat{\mu}_{lm}^{H} = 2\alpha \cdot \langle\psi_l, \psi_m {\rangle}_{H(\hat{K})},
\end{align}
where $\psi_l$ with $l=1,...,p$ denote the basis functions associated with
the reference element $\hat{K}$. Since the inner products are of integral type, the transformation theorem applies, analogously to \eqref{element_transformation_1}-\eqref{element_transformation_4}, for the mapping from the reference element to the local element, i.e$.$
\begin{align}
\mu_{lm}^{H,K} &= 2\alpha \cdot \langle \phi_l, \phi_m {\rangle}_{H(K)} \\
&= 2\alpha \cdot h_1^K \cdot ... \cdot h_n^K \langle\psi_l, \psi_m {\rangle}_{H(\hat{K})} \\
&= h_1^K \cdot ... \cdot h_n^K \hat{\mu}_{lm}^{H}.
\end{align}
Since the remaining derivation does not change, it is easy to see that formula \eqref{final_iters_estimate} is valid for the finite element discretised version of \eqref{abstract_min_problem}, when $f(\vec{u}_0)$ is replaced by $f(\vec{u}_0)-f(\vec{u}^*)$, and ${\lambda_{\max}^{\hat{M}}}, {\lambda_{\min}^{\hat{M}}}$ are replaced by ${\lambda_{\max}^{\hat{M}_H}}, {\lambda_{\min}^{\hat{M}_H}}$, respectively.


\section{Numerical experiments}\label{sec:experiments}

In this section, the results above are verified by solving the previously
discussed generic optimisation problem \eqref{abstract_min_problem} numerically. The behaviour predicted
is also observed to carry over to the numerical solution of a PDE constrained
optimisation problem in two dimensions. The simulations were performed using
the FEniCS automated finite element system
\cite{logg2012automated}. For the PDE constrained problem, the optimisations
were conducted according to the framework presented in \cite{funke2013framework}.
It makes use of the reduced problem formulation and of the automated optimisation
facilities of dolfin-adjoint \cite{Farrell2013}. The actual code employed to
conduct the simulations is in the supplement to this paper.


\subsection{Generic optimisation problem}

\begin{figure}[h]
    \centering
    \begin{subfigure}[b]{0.48\textwidth}
        \includegraphics[width=\textwidth]{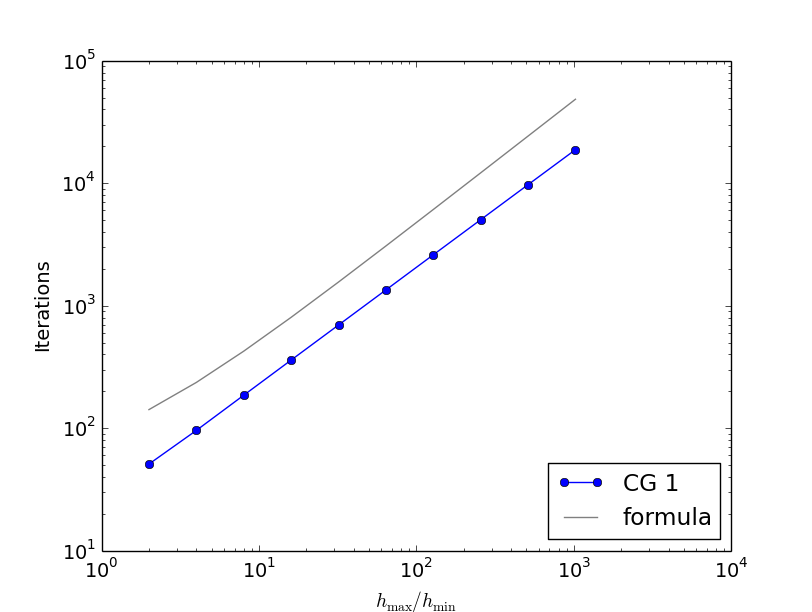}
        \caption{1-D}
        \label{1D}
    \end{subfigure} \\
    \begin{subfigure}[b]{0.48\textwidth}
        \includegraphics[width=\textwidth]{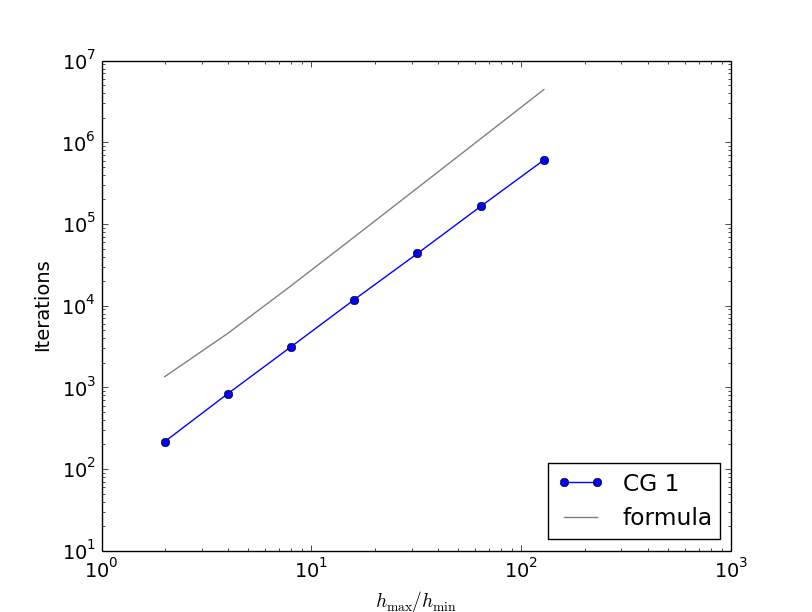}
        \caption{2-D}
        \label{2D}
    \end{subfigure}
    ~ 
    \begin{subfigure}[b]{0.48\textwidth}
        \includegraphics[width=\textwidth]{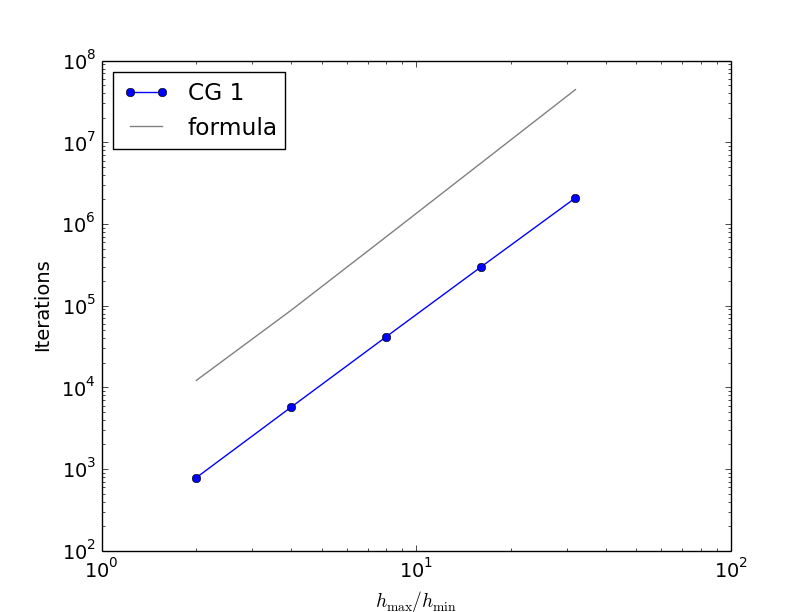}
        \caption{3-D}
        \label{3D}
    \end{subfigure}
    \vspace{0.4cm}
    \caption{Iterations of the steepest descent method with exact line search and $\ell^2$ represented gradients, using linear Lagrange finite elements on non-uniform refined meshes with respect to $h_{\max}/h_{\min}$ on a $\log$-$\log$ scale, and with error threshold $\varepsilon = 10^{-15}$ in the functional.}
    \label{loglog_plot}
\end{figure}

The iteration number estimate from the last section applies if the $\ell^2$ inner product is used in the representation of the derivative of the objective functional $f$. Corresponding simulations for the solution of \eqref{example_opti} with initial $u_0=0$ were run on successively refined meshes. The $(r+1)$-times refined mesh is derived from the $r$-times refined mesh by doubling the mesh resolution in part of the $r$-times refined mesh region. That way, $h_{\max}/h_{\min}$ doubles after every refinement. Figure \eqref{loglog_plot} displays the evolution of iteration numbers with respect to $h_{\max} / h_{\min}$ using linear Lagrange elements in one, two and three dimensions and the corresponding analytical iteration count estimate. The polynomial dependency between iteration number and $h_{\max}/h_{\min}$ in the estimate \eqref{final_iters_estimate} is numerically confirmed. Analogous dependencies were found for other and higher order finite elements implemented in FEniCS, emphasising that the iteration number evolution indeed exhibits the polynomial relationship stated in \eqref{final_iters_estimate}, independently of the finite element choice. \\

The iteration numbers from Figure \eqref{loglog_plot} compare to an iteration number that is always, and independently of the mesh or the dimension of the underlying domain, equal to $1$ when the functional derivative of $f$ is represented with respect to $L^2$ corresponding to the control variable space.\\

By \emph{uniform refinement} we mean that the mesh is refined by an equal factor throughout the domain such that the ratio $h_{\max}/h_{\min}$ remains unaltered. Consequently, we do not expect the iteration number to increase under uniform refinement according to \eqref{final_iters_estimate}. Indeed, the results given in table \ref{tab:iters_uniform_ref} confirm this. In fact, the iterations tend to decrease slightly for increasing refinement levels. This makes sense as the relative impact of the numerical error due to the boundaries gets smaller as the total number of elements increases. For the simulations, an initially non-uniform mesh was considered. Similar results are obtained for an initially uniform mesh, with overall smaller iteration numbers compared to table \ref{tab:iters_uniform_ref}. In summary, the increase of iteration counts depends on how the mesh is refined, and not on what initial mesh is used. Further, the iteration numbers increase for increasing orders of the Lagrange elements, since the larger the Lagrange finite elements order is, the larger is the ratio $h_{\max}/h_{\min}$.



\begin{table}[h]
\centering
\caption{Iterations of the steepest descent method with exact line search and $\ell^2$ represented gradients for Lagrange finite elements of different order (CG$1$-$5$) on uniformly refined meshes in $1$D depending on the refinement level (number of elements). The error threshold is given by $\varepsilon = 10^{-6}$ in the functional.}

\begin{tabular}{ @{} l @{}  @{} *5c @{}}    \\ \bottomrule
 \emph{Refinement \hspace{2mm}} & \emph{CG1} & \emph{CG2} & \emph{CG3} & \emph{CG4} & \emph{CG5} \\\midrule
 1 (96)           & 13        & 20        & 19        & 37        & 43 \\
 2 (192)          & 12        & 20        & 16        & 36        & 41        \\
 3 (384)          & 11        & 19        & 15        & 35        & 40        \\
 4 (768)          & 10        & 19        & 13        & 35        & 38        \\
 5 (1536)         & 8         & 19        & 13        & 35        & 37        \\ \bottomrule
\end{tabular}
\label{tab:iters_uniform_ref}
\end{table}

\subsection{Optimal control of the Poisson equation}


A common case for function based optimisation to occur is PDE-constraint optimisation. To exemplify this we consider a simple and generic problem in PDE-constraint optimisation, namely an optimal control problem constrained by the Poisson equation. Physically, this problem can the interpreted as finding the best heating/cooling of a hot plate to achieve a desired temperature profile. For this problem, we analyse the performance of multiple commonly used implementations of optimisation algorithms, based on different inner products. Mathematically, the problem is to minimise the following tracking type functional
\begin{align}
\min_{f \in H} \hspace{1mm} \frac{1}{2} \|u - d\|_{L^2(\Omega)}^2 + \frac{\alpha}{2} \| f \|_{H}^2,
\label{poisson-mother-problem1}
\end{align}
subject to the Poisson equation with Dirichlet boundary conditions,
\begin{align}
\Delta u &= f  \text{ \hspace{2mm}in } \Omega \label{poisson-mother-problem2}\\
u &= 0  \text{ \hspace{2mm}on } \partial\Omega,
\label{poisson-mother-problem3}
\end{align}
where $\Omega = [0,1]^2$, $H=L^2(\Omega)$ or $H^1(\Omega)$,  $u: \Omega \to \mathbb{R}$ is the unkown temperature, $f \in H$ is the unknown control function acting as source term ($f(x) > 0$ corresponds to heating and $f(x) < 0$ corresponds to cooling), $d \in L^2(\Omega)$ is the given desired temperature profile and $\alpha \in [0, \infty)$ is a Tikhonov regularisation parameter. For a proof of existence and uniqueness of this optimisation problem, we refer to section 1.5 of \cite{pinnau2008optimization}.\\

All implementations apply the limited memory BFGS method in order to approximate Hessians using information from up to a certain number of previous iterations. The L-BFGS-B algorithm, derived in \cite{byrd1995limited} and implemented as part of Scipy's optimisation package \cite{scipy2001}, is a limited memory BFGS algorithm with support for bounds on the controls. The Limited Memory, Variable Metric (LMVM) method from the TAO library \cite{benson2005tao} uses a BFGS approximation with a Mor\'e-Thuente line search. The interior point method described in \cite{wachter2006implementation} and implemented in IPOPT \cite{kawajir2015} is an interior point line search filter method. Here, the BFGS method approximates the Hessian of the Lagrangian and not of the functional itself. \\

These optimisation method implementations are insensitive towards changes in the inner product space, but instead assume an Euclidean space. This compares to the BFGS implementation from Moola (\cite{funke2016}), that allows the user to specify the inner products. By applying Riesz maps associated to the underlying Hilbert space, the Moola respects inner products.\\

In table \eqref{pde-constraint_table}, iteration counts for the different algorithm implementations solving \eqref{poisson-mother-problem1} subject to \eqref{poisson-mother-problem2} and \eqref{poisson-mother-problem3} are displayed for $H=L^2$ and $H=H^1$, respectively. In both cases, the inner product insensitive methods exhibit mesh dependence. This can be seen from the increase in iteration numbers for increasing ratio between largest and smallest element size $h_{\max}/h_{\min}$, reflecting the non-uniformity in the discretisation. Iteration counts increase roughly by factors between $4$-$25$ when comparing the least and most non-uniform meshes. Relating this to large scale optimisation problems, it is obvious that the ability of efficiently dealing with non-uniformity in the discretisation can very well decide upon the computational feasibility of solving the problem. The optimisation method that respects the inner product of the control space is approximately mesh independent.





\begin{table}[h]
\centering
\caption{Iterations of different common optimisation method implementations in Python using different inner products. The error threshold is given by $\varepsilon = 10^{-7}$ in the gradient norm $\|\mathcal{R}_H(f')\|_H$. History length for Hessian approximation equal to 10, except for LMVM in TAO, where it is equal to 5. \\}
\begin{subtable}{1.\textwidth}
\centering
\begin{tabular}{ @{} c @{} l | l @{} @{} *6c @{}}     \bottomrule
 inner product \hspace{2mm} & 	Implementation	&\emph{$h_{\max}/h_{\min}$ \hspace{9mm}} & 4 & 8 & 16 & 32 & 64 	&	128		 \\ \midrule
$\ell^2$ & 	Scipy	&  BFGS            & 27        & 47        & 75        & 97        & 133 	 &	189 \\
$\ell^2$ & 	TAO	& LMVM           	& 42        & 60        & 77        & 111        & 131    &   155 \\
$\ell^2$ & 	IPOPT	& Interior Point           & 28        & 38        & 65        & 88        & 117    & 139   \\ \midrule
$L^2$ & 	Moola		& BFGS         & 22        & 20        & 22       & 23        & 23    &   	27 \\  \bottomrule
\end{tabular}
 \caption{$H=L^2(\Omega)$}
\end{subtable}

\begin{subtable}{1.\textwidth}
\centering
\begin{tabular}{ @{} c @{} l | l @{} @{} *6c @{}}     \bottomrule
 inner product \hspace{2mm} & Implementation	&\emph{$h_{\max}/h_{\min}$ \hspace{9mm}} & 4 & 8 & 16 & 32 & 64 	&	128		 \\ \midrule
$\ell^2$ & 	Scipy		& BFGS          & 25        & 53        & 118        & 222        & 244 	 &	594 \\
$\ell^2$ & 	TAO		& LMVM          	& 38        & 38        & 110        & 144        & 141    &   384 \\
$\ell^2$ & 	IPOPT		& Interior Point          & 31        & 38        & 85        & 127        & 182    & 443   \\ \midrule
$H^1$ & 		Moola		& BFGS         & 4        & 4        & 4        & 4        & 4    &   	4 \\ \bottomrule
\end{tabular}
 \caption{$H=H^1(\Omega)$}
\end{subtable}
\label{pde-constraint_table}
\end{table}

\section{Conclusion}

This work can be seen as a first step in quantifying the computational expense of disrespecting inner products. It illustrates the importance of employing the inner product induced by the control space for optimisations by deriving a rigorous estimate of the cost of failing to do so for a simple and generic problem. Simulations confirm the derived dependency of iteration counts mesh non-uniformity, while similar numerical results are found for a more complex PDE constrained optimisation problem. These results provide an intuitive and analytically coherent model of mesh dependence and its corresponding errors for this class of problems. The ability to employ mesh refinement to resolve features of interest is a core advantage of the finite element method. Effectively dealing with non-uniformity in the mesh should therefore be a core capability of the optimisation method.



\section*{Acknowledgements}
The authors would like to thank Colin J$.$ Cotter and Magne Nordaas for valuable discussions and collaboration. This work was supported by the Grantham Institute for Climate Change, the Engineering and Physical Sciences Research Council [Mathematics of Planet Earth CDT, grant number EP/L016613/1], the Natural Environment Research Council [grant number NE/K008951/1], and the Research Council of Norway through a Centres of Excellence grant to the Center for Biomedical Computing at Simula Research Laboratory (project number 179578).

\nocite{*}
\bibliographystyle{siamplain}
\bibliography{schwedes_etal}


\end{document}